%% file: weakly.tex
\documentclass[a4paper]{amsart}

\usepackage{amssymb,pstricks,amscd}
\usepackage[totalwidth=17cm,totalheight=24cm]{geometry}
\usepackage{graphicx}
\input psfig

\begin{document}

\input{macros}

\title{On the infinitesimal rigidity of polyhedra with vertices in 
convex position}
\author{Ivan Izmestiev}
\thanks{The first-named author was supported by the DFG Research Unit 
``Polyhedral Surfaces''.}
\address{Institut f\"ur Mathematik, MA 8-3\\
Technische Universit\"at Berlin\\
Strasse des 17. Juni 136\\
10623 Berlin, Germany}
\email{izmestiev@math.tu-Berlin.de}
\author{Jean-Marc Schlenker}
\address{Institut de Math\'ematiques de Toulouse, UMR CNRS 5219\\
Universit\'e Toulouse III\\
31062 Toulouse cedex 9, France}
\email{schlenker@math.ups-tlse.fr}
\date{October 2007 (v1)}

\begin{abstract}
Let $P \subset \R^3$ be a polyhedron. It was conjectured
that if $P$ is weakly convex (i.~e. its vertices lie on the boundary of 
a strictly convex domain) and decomposable
(i.~e. $P$ can be triangulated without adding new vertices), then it is
infinitesimally rigid. We prove this conjecture under a weak additionnal assumption
of codecomposability.

The proof relies on a result of independent
interest concerning the Hilbert-Einstein function of a triangulated convex polyhedron. We determine the signature of the Hessian of that function with respect to deformations of the interior edges. In particular, if there are no interior vertices, then the Hessian is negative definite.
\end{abstract}
\maketitle


\section{Introduction}

\subsection{The rigidity of convex polyhedra.}

The rigidity of convex polyhedra is a classical result in geometry,
first proved by Cauchy \cite{cauchy} using ideas going back to Legendre \cite{legendre}.

\begin{thm}[Cauchy 1813, Legendre 1793] \label{thm:Cauchy-Legendre}
Let $P,Q\subset \R^3$ be two convex polyhedra with the same combinatorics and such
that corresponding faces are isometric. Then $P$ and $Q$ are congruent.
\end{thm}

This result had a profound influence on geometry over the last two centuries, it
led for instance to the discovery of the rigidity of smooth convex surfaces in
$\R^3$, to Alexandrov's rigidity and to his results on the realization of
positively curved cone-metrics on the boundary of polyhedra (see \cite{alex}).

From a practical viewpoint, global rigidity is perhaps not as relevant as 
infinitesimal rigidity (see Section \ref{subsec:def} for a definition). Although the infinitesimal rigidity of convex polyhedra can be proved using Cauchy's argument,
the first proof was given much later by Dehn \cite{dehn}, and is completely 
different from Cauchy's.

\begin{thm}[Dehn 1916] \label{tm:dehn}
Any convex polyhedron is infinitesimally rigid.
\end{thm}

In this paper we deal with a generalization of Theorem \ref{tm:dehn} to a 
vast class of non-convex polyhedra. The main idea is that it is not
necessary to consider convex polyhedra, what is important is that the
vertices are in convex position. Additional assumptions are necessary,
however, that are automatically satisfied for convex polyhedra.

\subsection{Main result.}
By a polyhedron we mean a body in $\R^3$ bounded by a closed polyhedral
surface.

\begin{figure}[ht]
\includegraphics[width=0.8\textwidth]{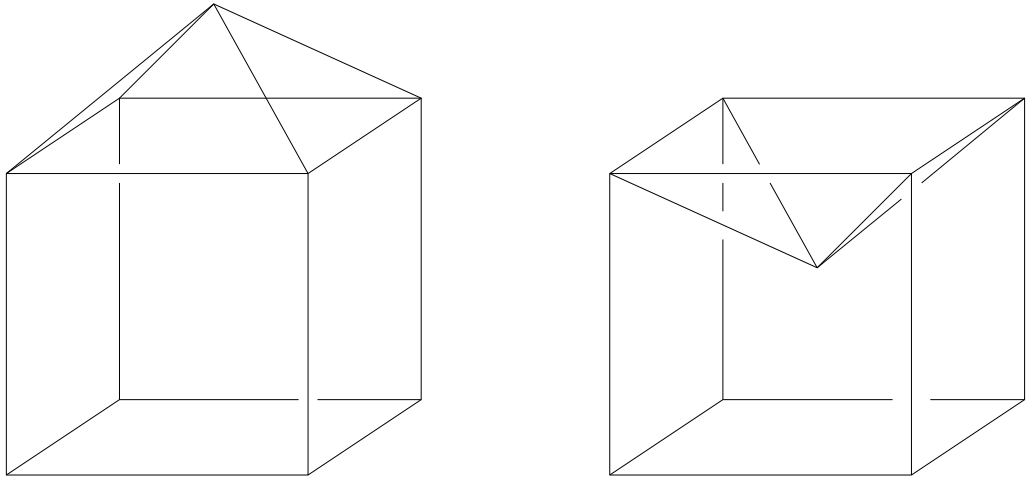}
\caption{Cauchy's rigidity theorem fails for non-convex polyhedra.}
\label{fig:CL-counter}
\end{figure}

\begin{figure}[ht]
\includegraphics[height=0.35\textwidth]{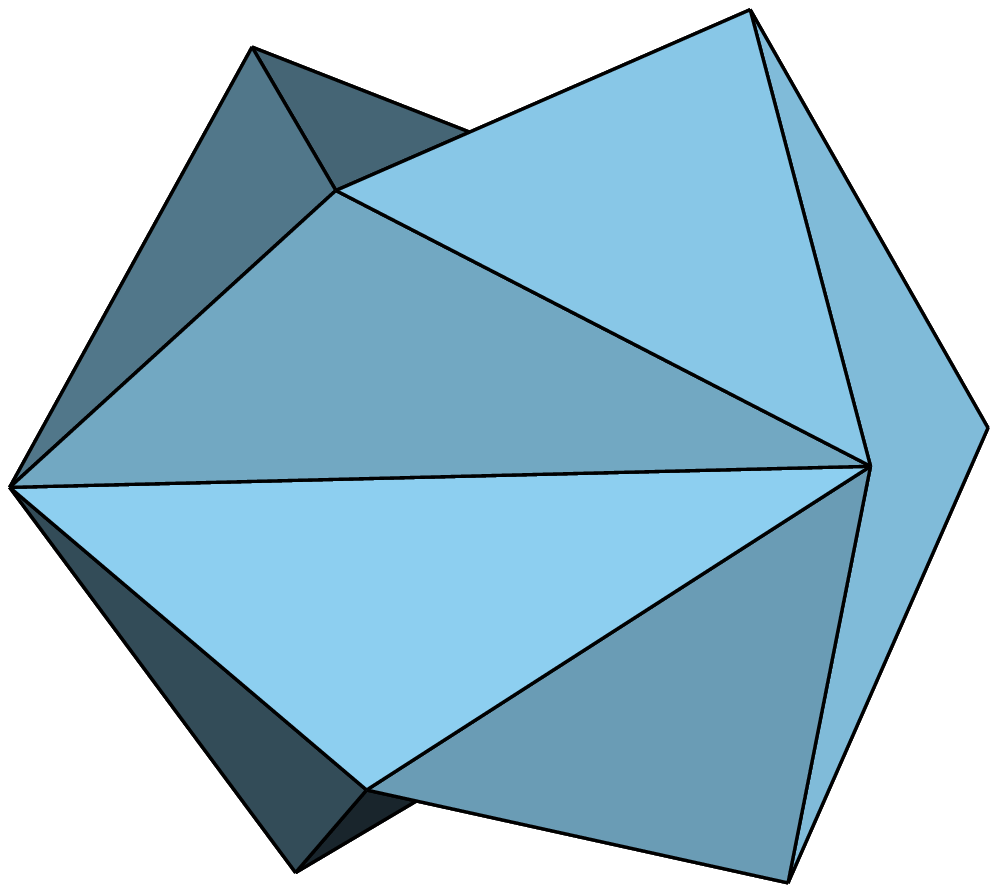} \hspace{0.1\textwidth}
\includegraphics[height=0.35\textwidth]{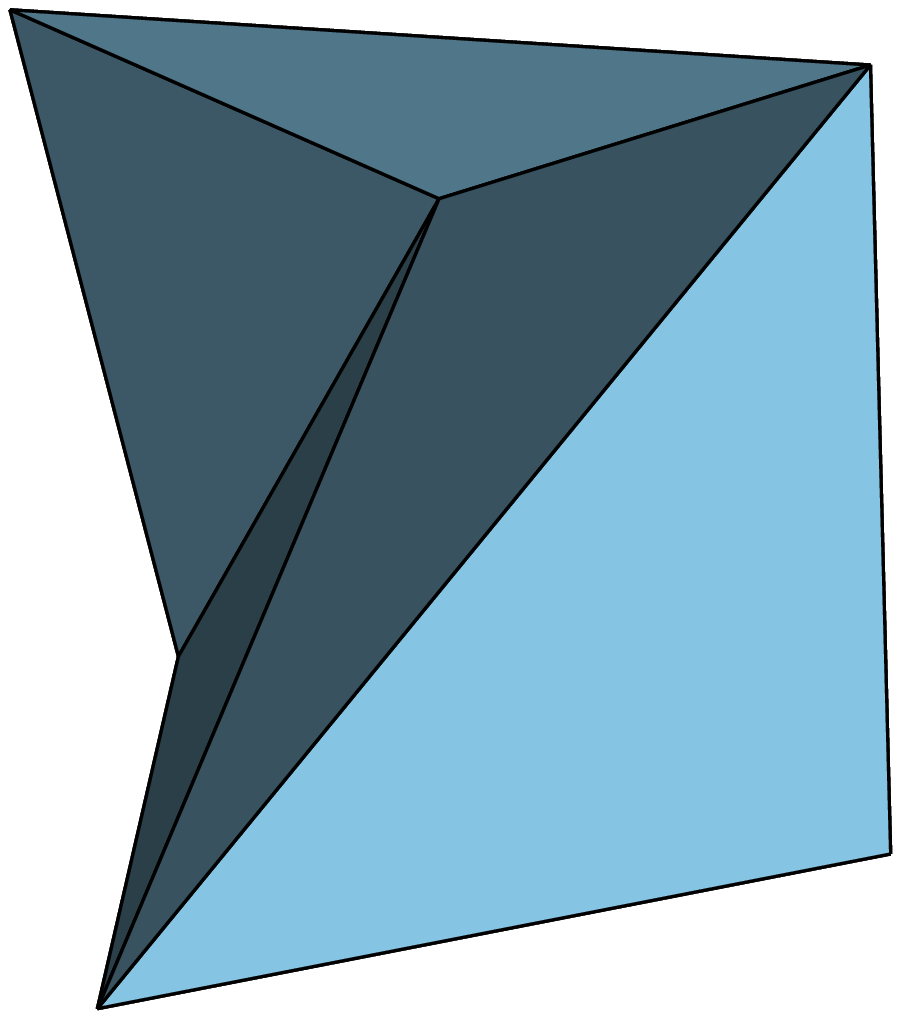}
\caption{Examples of infinitesimally flexible polyhedra: Jessen's orthogonal icosahedron \cite{Jes67}, \cite{MW07} and Wunderlich's twisted octahedron \cite{Wun65}. Both are weakly convex (Definition \ref{dfn:wc}) but not decomposable (Definition \ref{dfn:dec}).}
\label{fig:IcoOcta}
\end{figure}

For non-convex polyhedra neither of the theorems \ref{thm:Cauchy-Legendre} and \ref{tm:dehn} is true. It is easy to find a counterexample to Theorem \ref{thm:Cauchy-Legendre}, see Figure \ref{fig:CL-counter}. Counterexamples to Theorem \ref{tm:dehn} are more complicated, see Figure \ref{fig:IcoOcta}.

\begin{defi}
\label{dfn:wc}
A polyhedron $P\subset \R^3$ is called \emph{weakly convex} if its
vertices are in convex position in $\R^3$.
\end{defi}

In other words, $P$ is weakly convex if its vertices are the vertices of a strictly convex polyhedron.

\begin{defi}
\label{dfn:dec}
A polyhedron $P$ is called \emph{decomposable} if it can be triangulated without adding new vertices.
\end{defi}

In other words, every simplex of the triangulation must have vertices among those of $P$.

Our work was motivated by the following conjecture.

\begin{conj} \label{cj:main}
Every weakly convex decomposable polyhedron is infinitesimally rigid.
\end{conj}

Any infinitesimally flexible polyhedron known to us fails to satisfy one of the assumptions of Conjecture~\ref{cj:main}. Thus, both polyhedra on Figure \ref{fig:IcoOcta} are weakly convex but not decomposable. An infinitesimally flexible non-convex octahedron pictured in \cite{Glu75} is decomposable but not weakly convex.

The main result of this paper is the proof of a weakening of Conjecture \ref{cj:main}.

\begin{thm} \label{tm:main}
Let $P$ be a weakly convex, decomposable, and weakly codecomposable 
polyhedron. Then $P$ is infinitesimally rigid.
\end{thm}

Note that $P$ is not required to be homeomorphic to a ball. The hypothesis that 
$P$ is weakly codecomposable, however, appears to be quite weak 
for polyhedra homeomorphic to a ball. In the appendix we describe a simple
example of a polyhedron which is not weakly codecomposable, however it is
not homeomorphic to a ball; it's quite possible that this example can be modified
fairly simply to make it contractible.

\begin{defi}
\label{dfn:codec}
A polyhedron $P$ is called {\it codecomposable} if its complement in $\conv P$ 
can be triangulated without adding new vertices. $P$ is called 
{\it weakly codecomposable} if $P$ is contained in a convex polyhedron $Q$,
such that all vertices of $P$ are vertices of $Q$ and that the complement
of $P$ in $Q$ can be triangulated without adding new vertices.
\end{defi}

The examples of non-codecomposable weakly convex polyhedra homeomorphic to a
ball are quite complicated 
\cite{aichholzer}, so a counterexample to Conjecture \ref{cj:main} would be 
difficult to construct. On the other hand, the codecomposability assumption 
is used in our proof of Theorem \ref{tm:main} in a very essential way. Thus the 
question whether the codecomposability assumption may be omitted remains wide open.

\subsection{Earlier results.}

Conjecture \ref{cj:main} originated as a question in \cite{rcnp}, where a related result
was proved: if $P$ is a decomposable polyhedron such that there exists an 
ellipsoid which intersects all edges of $P$ but contains none of its 
vertices, then $P$ is infinitesimally rigid. The proof relies on hyperbolic geometry,
more precisely the properties of the volume of hyperideal hyperbolic polyhedra.

Two special cases of the conjecture were then proved in \cite{vienna}: when $P$
is a weakly convex suspension containing its north-south axis, and when $P$ has only one concave
edge, or two concave edges adjacent to a vertex. The proof for suspensions was
based on stress arguments, while the proof of the other result used a refinement
of Cauchy's argument.

More recently, it was proved in \cite{star} that the conjecture holds when $P$ is
star-shaped with respect to one of its vertices. This implies the two results in
\cite{vienna}. The proof was based on recent results of \cite{izmestiev} concerning
convex caps. The result of \cite{star} -- and therefore the two results of 
\cite{vienna} -- are consequences of Theorem \ref{tm:main}, since it is not difficult
to show that a polyhedron which is star-shaped with respect to one of its 
vertices is codecomposable (the proof actually appears as a step in \cite{star}).

\subsection{Definitions.}
\label{subsec:def}
Let us set up the terminology.

Every polyhedron has faces, edges, and vertices. (Note that a face can be multiply connected.) For every polyhedron $P$ its boundary $\partial P$ can be triangulated without adding new vertices. This follows from the well-known fact that every polygon, be it non-convex or multiply connected, can be cut into triangles by diagonals.

\begin{defi}
\label{dfn:InfRig}
Let $P$ be a polyhedron, and let $S$ be a triangulation of its boundary with $\V(S)=\V(P)$. The polyhedron $P$ is called \emph{infinitesimally rigid} if every infinitesimal flex of the 1-skeleton of $S$ is trivial.
\end{defi}

By an \emph{infinitesimal flex} of a graph in $\R^3$ we mean an assignment of vectors to the vertices of the graph such that the displacements of the vertices in the assigned directions induce a zero first-order change of the edge lengths:
$$
(p_i - p_j) \cdot (q_i - q_j) = 0 \mbox{ for every edge }p_ip_j,
$$
where $q_i$ is the vector associated to the vertex $p_i$. An infinitesimal flex is called trivial if it is the restriction of an infinitesimal rigid motion of $\R^3$.

It is not hard to see that for any two triangulations $S$ and $S'$ as in Definition \ref{dfn:InfRig} a flex of $S$ is a flex of $S'$. Thus the infinitesimal rigidity of a polyhedron is a well-defined notion.

\begin{figure}[ht]
\includegraphics[width=0.3\textwidth]{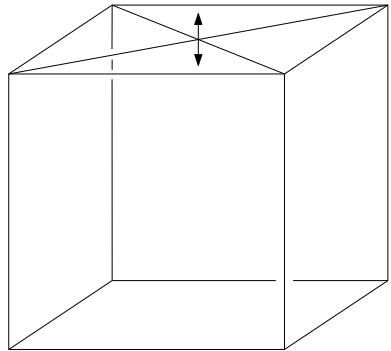}
\caption{Convex triangulated sphere with a flat vertex. Moving this vertex in the vertical direction produces a non-trivial infinitesimal flex.}
\label{fig:FlatVert}
\end{figure}

We will need to deal with triangulations of the boundaries of polyhedra that contain additional vertices. Infinitesimal flexes and infinitesimal rigidity for such \emph{triangulated spheres} are defined in the same way. Note that a triangulated sphere may be flexible even if it bounds a convex polyhedron, see Figure \ref{fig:FlatVert}.

\begin{defi}
Let $S$ be a triangulation of the boundary of a polyhedron $P$. A vertex $p$ of $S$ is called a \emph{flat vertex} if it lies in the interior of a face of $P$.
\end{defi}

The following statement, which generalizes Dehn's theorem to convex simplicial spheres, 
will be useful below. Its proof follows directly from Dehn's Theorem, by considering
the restriction of an infinitesimal flex to the non-flat vertices.

\begin{thm}
\label{thm:DehnGen}
Let $S$ be a triangulation of the boundary of a convex polyhedron $P$. Then every infinitesimal flex of $S$ is the sum of an infinitesimal rigid motion and of displacements of flat vertices in the directions orthogonal to their ambient faces.
\end{thm}

\subsection{The Hilbert-Einstein function.}
The proof of Theorem \ref{tm:main} is based on some striking properties of the discrete
Hilbert-Einstein function, also known in the physics community as the Regge 
function \cite{regge}. First we have to define a space of deformations
of a triangulated polyhedron.

\begin{defi}
Let $T$ be a triangulation of a polyhedron $P$, and let $e_1,\cdots, e_n$ be the interior edges of $T$. We denote by $\cD_{P,T}$ the space of $n$-tuples $(l_1, \cdots, l_n) \in \R_{>0}^n$ such that for every simplex $\sigma$ of $T$, replacing the lengths of the edges of $\sigma$ that are interior edges of $T$ by the corresponding $l_j$ produces a non-degenerate simplex. 
\end{defi}

For every element $l\in \cD_{P,T}$ there is an associated metric structure on $P$ obtained by gluing the simplices with changed edge lengths. The resulting metric space is locally Euclidean except that it has cone singularities
along the interior edges of $T$. For every $i\in \{ 1,\cdots, n\}$, denote by $\omega_i$ the total angle around $e_i$ and by $\kappa_i :=2\pi-\omega_i$ the singular curvature along $e_i$. Let $e'_1, \cdots, e'_r$ be the boundary
edges of $P$; for every $j\in \{ 1,\cdots, r\}$ denote by $\alpha_j$ the dihedral angle of $P$ at $e'_j$, and by $l'_j$ the length of $e'_j$. 

\begin{defi}
The \emph{Hilbert-Einstein function} on $\cD_{P,T}$ is given by the formula
$$ \cS(l) := \sum_{i=1}^n l_i \kappa_i + \sum_{j=1}^r l'_j (\pi-\alpha_j)~. $$  
\end{defi}

\subsection{The Schl\"afli formula.}
It is a key tool in polyhedral geometry. It has several generalizations, but
the 3-dimensional Euclidean version states simply that, under a first-order
deformation of any Euclidean polyhedron, 
\begin{equation}
  \label{eq:schlafli}
  \sum_e l_e d\alpha_e = 0~,
\end{equation}
where the sum is taken over all edges $e$, with $l_e$ denoting the length of the edge $e$, and 
$\alpha_e$ the dihedral angle at $e$. This equality is also known as the
``Regge formula''.

It follows directly from the Schl\"afli formula that, under any first-order
variation of the lengths of the interior edges of a triangulation $T$ of
the polyhedron $P$ --- i.e., for any tangent vector to $\cD_{P,T}$ --- the first-order
variation of $\cS$ is simply
\begin{equation}
  \label{eq:cor-schlafli}
  d\cS = \sum_{i=1}^n \kappa_i dl_i~.
\end{equation}

As a consequence, the Hessian of $\cS$ equals the Jacobian of the map $(l_i)_{i=1}^n \mapsto (\kappa_i)_{i=1}^n$.

\begin{defi}
\label{dfn:M_T}
Let $T$ be a triangulation of a polyhedron $P$ with $n$ interior edges. Define the $n \times n$ matrix $M_T$ as
$$
M_T = \left( \frac{\partial \omega_i}{\partial l_j} \right) = -\left( \frac{\partial^2 \cS}{\partial l_i \partial l_j} \right).
$$
The derivatives are taken at the point $l \in \cD_{P,T}$ that corresponds to the actual edge lengths in $T$.
\end{defi}

The arguments in this paper use only $M_T$, and not directly the Hilbert-Einstein function $\cS$. The fact that $M_T$ is minus the Hessian of $\cS$ does imply, however, that $M_T$ is symmetric.

The matrix $M_T$ is directly related to the
infinitesimal rigidity of $P$, an idea that goes back to Blaschke and Herglotz.
\footnote{Blaschke and Herglotz suggested that the critical points of the Hilbert-Einstein 
function on a manifold with boundary (in the smooth case), with fixed boundary metric, 
correspond to Einstein metrics,
i.e., to constant curvature metrics in dimension $3$. The analog of $M_T$ in this
context is the Hessian of the Hilbert-Einstein function.}

\begin{lemma} \label{lm:nondegen}
Let $T$ be a triangulation of a polyhedron $P$ without interior vertices.
Then $P$ is infinitesimally rigid if and only if $M_T$ is non-degenerate.
\end{lemma}

The proof can be found in \cite{bobenko-izmestiev,star}.
It is based on the observation that an isometric deformation of $P$
induces a first-order variation of the interior edge lengths but a zero variation of the angles around them.

\subsection{The second-order behavior of $\cS$.}
The following is the key technical statement of the paper.

\begin{thm} \label{tm:novertex}
Let $P$ be a convex polyhedron, and let $T$ be a triangulation of $P$ with $\V(T) = \V(P)$. Then $M_T$ is positive definite.
\end{thm}

Theorem \ref{tm:novertex} is actually a special case of the following theorem that describes
the signature of $M_T$ for $T$ \emph{any} triangulation of $P$. 

\begin{thm} \label{tm:withvertices}
Let $P$ be a convex polyhedron, and let $T$ be a triangulation of $P$
with $m$ interior and $k$ flat vertices. Then the dimension of the kernel of $M_T$
is $3m + k$, and $M_T$ has $m$ negative eigenvalues.
\end{thm}

\subsection{From Theorem \ref{tm:novertex} to Theorem \ref{tm:main}.}

Let $P$ be a polyhedron that satisfies the assumptions of Theorem~\ref{tm:main}. 
Since $P$  is decomposable and weakly codecomposable, there exists a convex 
polyhedron $Q$ such that all vertices of $P$ are vertices of $Q$, and a 
a triangulation $\overline{T}$ of $Q$ that contains a triangulation $T$ of $P$ 
and whose vertices are only the vertices of $Q$. 
It is easy to see that the matrix $M_T$ is then a principal minor of the matrix $M_{\overline{T}}$. 
By Theorem \ref{tm:novertex}, $M_{\overline{T}}$ is positive definite, thus $M_{T}$ is 
positive definite as well. In particular, $M_T$ is non-degenerate. Lemma \ref{lm:nondegen} 
implies that the polyhedron $P$ is infinitesimally rigid.

Since Theorem \ref{tm:novertex} is a special case of Theorem \ref{tm:withvertices}, 
the rest of this paper deals with proving Theorem \ref{tm:withvertices}.

\subsection{Proving Theorem \ref{tm:withvertices}}
\label{subsec:plan}

The proof is based on a standard procedure. In order to show that the matrix $M_T$ has the desired property for every triangulation $T$, we prove three points:
\begin{itemize}
 \item any two triangulations can be connected by a sequence of moves;
 \item the moves don't affect the desired property;
 \item the property holds for a special triangulation.
\end{itemize}
These points are dealt with in the given order in the next three sections of the paper.

\section{Connectedness of the set of triangulations}
\label{sec:conn}
Moves on simplicial complexes are well-studies, see an overview in \cite{Lick99}. There are several theorems stating that any two triangulations of a given manifold can be connected by certain kinds of simplicial moves. Note, however, that we are in a different situation here, since we deal with triangulations of a fixed geometric object. Taking a closer look, one sees that a simplicial move is defined as a geometric move preceded and followed by a simplicial isomorphism. Performing an isomorphism is the possibility that is missing in our case.

To emphasize the difference between the combinatorial and the geometric situation, let us cite a negative result concerning geometric moves. Santos \cite{San05} exhibited two triangulations with the same set of vertices in $\R^5$ that cannot be connected via $2 \leftrightarrow 5$ and $3 \leftrightarrow 4$ bistellar moves. For an overview on geometric bistellar moves see \cite{San06}.

\subsection{Geometric stellar moves: the Morelli-W{\l}odarczyk theorem}
A positive result on geometric simplicial moves was obtained by Morelli \cite{morelli} and W{\l}odarczyk \cite{wlodarczyk}. As a crucial step in the proof of the weak Oda conjecture, they showed that any two triangulations of a convex polyhedron can be connected by a sequence of \emph{geometric stellar moves}.

\begin{figure}[ht]
\includegraphics[width=0.8\textwidth]{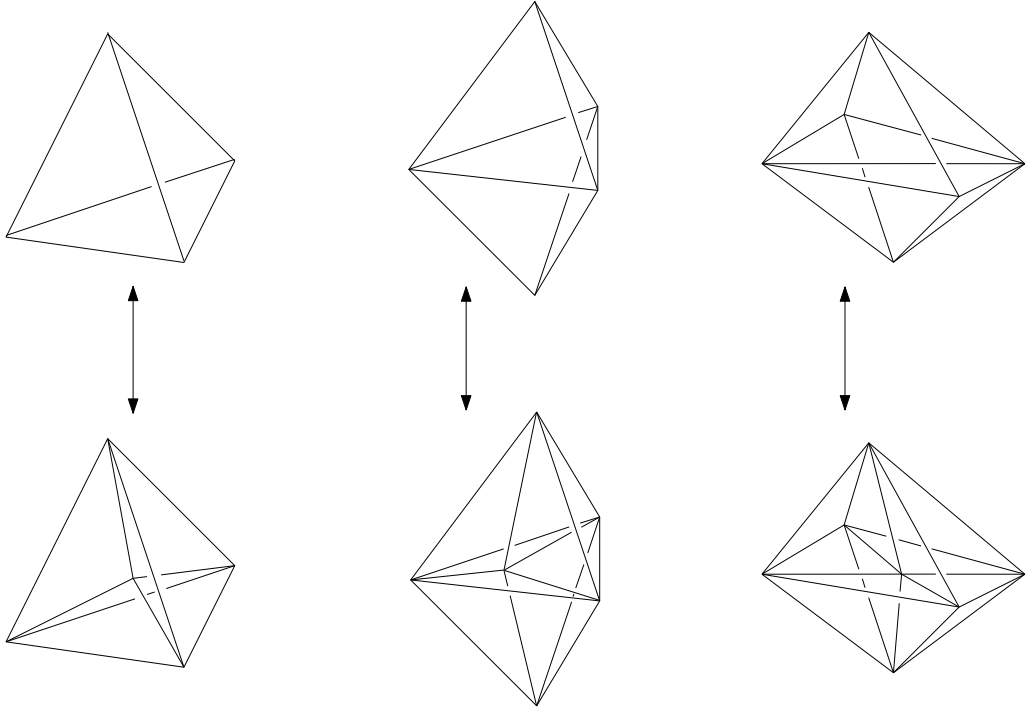}
\caption{Interior stellar moves in dimension $3$.}
\label{fig:StarWeld}
\end{figure}

\begin{figure}[ht]
\includegraphics[width=0.6\textwidth]{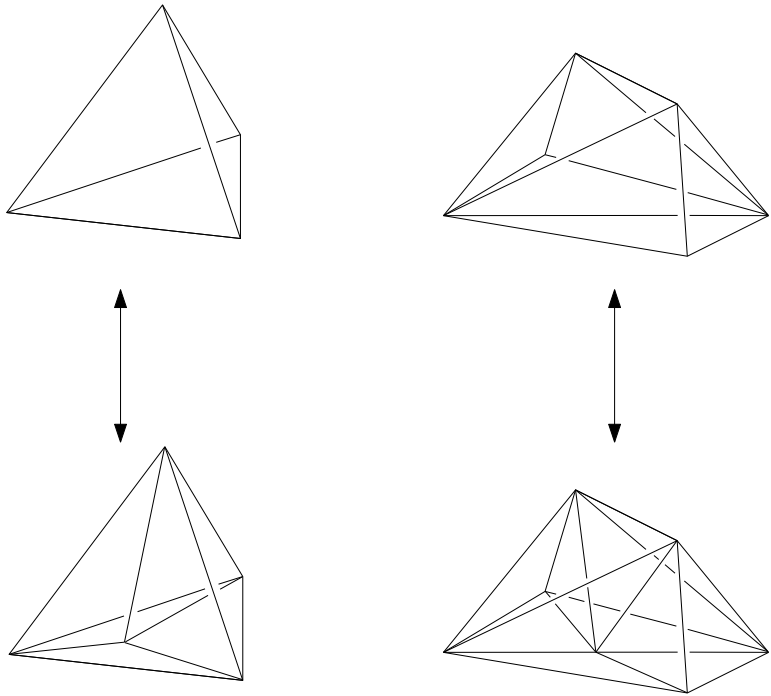}
\caption{Boundary stellar moves in dimension $3$.}
\label{fig:dStarWeld}
\end{figure}

\begin{defi}
Let $p$ be an interior point of a simplex $\sigma \subset \R^n$. The \emph{starring} of $\sigma$ at $p$ is an operation that replaces $\sigma$ by the cone with the apex $p$ over the boundary of $\sigma$.

Let $T$ be a triangulation of a subset of $\R^n$, let $\sigma$ be a simplex of $T$, and let $p$ be a point in the relative interior of $\sigma$. The operation of \emph{starring} of $T$ at $p$ consists in replacing the star $\st \sigma$ of $\sigma$ by the cone with apex $p$ over the boundary of $\st \sigma$. The operation inverse to starring is called \emph{welding}.

Starrings and weldings are called \emph{stellar moves}.
\end{defi}
See Figures \ref{fig:StarWeld} and \ref{fig:dStarWeld} for stellar moves in dimension $3$. Figure \ref{fig:StarWeld} depicts starrings and weldings at interior points of $T$, Figure \ref{fig:dStarWeld} --- starrings and weldings at boundary points. Note that in the case of a boundary point our definition is not completely correct: a starring replaces $\st \sigma$ by the cone over $\partial\, \st \sigma \setminus \partial T$.

\begin{thm}[Morelli \cite{morelli}, W{\l}odarczyk \cite{wlodarczyk}]
\label{thm:MW}
Any two triangulations of a convex polyhedron $P \subset \R^n$ can be connected by a sequence of geometric stellar moves.
\end{thm}

We give here an outline of the Morelli's proof using more elementary language and tools.
\smallskip

\noindent\textit{Outline of the proof.} Let $T$ and $T'$ be two triangulations of $P$. A triangulation $\Sigma$ of $P \times [0,1]$ with $\Sigma_{P \times \{1\}} = T$ and $\Sigma|_{P \times \{0\}} = T'$ is called a \emph{simplicial cobordism} between $T$ and $T'$.

\begin{defi}
Let $\pr$ denote the orthogonal projection $P \times [0,1] \to P$. A simplex $\sigma \in \Sigma$ is called a \emph{circuit} if $\dim \pr(\sigma) < \dim \sigma$ and $\sigma$ is inclusion-minimal with this property.
\end{defi}

Clearly, the stars of the circuits are simplicial balls with no vertical faces and
$$
\Sigma = \bigcup_{\sigma \mbox{ \small{circuit}}} \st \sigma
$$
with disjoint interiors.

\begin{defi}
A simplicial cobordism $\Sigma$ is called \emph{collapsible} if there is a sequence of triangulations $\Sigma = \Sigma_0, \Sigma_1, \ldots, \Sigma_N = \Sigma_{P \times \{0\}}$ such that
\begin{itemize}
 \item $\Sigma_{i+1} = \Sigma_i \setminus \st \sigma_i$ for a circuit $\sigma_i$;
 \item the upper boundary of $\Sigma_i$ projects one-to-one on $P$ for every $i$.
\end{itemize}
\end{defi}
In other words, $\Sigma$ is collapsible if it can be ``dismantled with a crane''.

\begin{lemma}
The triangulation $\pr(\partial^+ \Sigma_{i+1})$ can be obtained from $\pr(\partial^+ \Sigma_i)$ by a starring with a subsequent welding. Here $\partial^+$ denotes the upper boundary.
\end{lemma}
\begin{proof}
For every circuit $\sigma$ the transformation $\pr(\partial^+ \sigma) \rightsquigarrow \pr(\partial^- \sigma)$ is a bistellar move and can be realized by a starring and a welding. These extend to a starring and a welding in $\st \sigma$.
\end{proof}

Thus, a collapsible simplicial cobordism between two triangulations gives rise to a sequence of stellar moves joining the triangulations.

\begin{defi}
A triangulation $\Sigma$ is called \emph{coherent} if there is a function $h: |\Sigma| \to \R$ that is piecewise linear with respect to $\Sigma$ and strictly convex across every facet of $\Sigma$. (Here $|\Sigma| = \cup_{\sigma \in \Sigma} \sigma$ is the support of $\Sigma$.)
\end{defi}

As an example, the reader can check that the barycentric subdivision of any
polytope $P$ is coherent -- one can choose $h$ taking on each vertex of the
barycentric subdivision an integer value equal to the dimension of the
corresponding face of $P$.

\begin{lemma}
A coherent simplicial cobordism is collapsible.
\end{lemma}

\begin{proof}
Let $\sigma$ and $\sigma'$ be two circuits of $\Sigma$ such that some point of $\st \sigma$ lies directly above a point of $\st \sigma'$. It follows that $\frac{\partial h}{\partial x_{n+1}}|_\sigma > \frac{\partial h}{\partial x_{n+1}}|_{\sigma'}$, where $\frac{\partial}{\partial x_{n+1}}$ denotes the derivative in the vertical direction. Thus the stars of the circuits can be lifted up in the non-decreasing order of the vertical derivative of $h$ on the circuits.
\end{proof}

To prove the theorem, we construct a coherent cobordism between stellar subdivisions of $T$ and $T'$. (By a stellar subdivision we mean the result of a sequence of starrings.)

\begin{lemma}[]
\label{lem:Subd}
Let $\Sigma$ and $\Sigma'$ be two triangulations with the same support. Then $\Sigma$ can be stellarly subdivided to a triangulation $\Sigma''$ that refines $\Sigma'$.
\end{lemma}

The reader can find a proof of this classical statement in e.g. \cite{glaser}.

\begin{lemma}
\label{lem:SSubd}
Let $\Sigma$ be an arbitrary triangulation. 
Then $\Sigma$ can be stellarly subdivided to a coherent triangulation.
\end{lemma}

\begin{proof}
By Lemma \ref{lem:Subd}, the barycentric triangulation of the support $|\Sigma|$ of $\Sigma$
can be stellarly subdivided to a triangulation $\Sigma'$ that refines $\Sigma$. 
But the barycentric subdivision of any polytope is coherent. 
Since starring a coherent triangulation produces a coherent triangulation, $\Sigma'$ is coherent.

Again by Lemma \ref{lem:Subd}, the triangulation $\Sigma$ can be stellarly subdivided to a 
triangulation $\Sigma''$ that refines $\Sigma'$. We claim that $\Sigma''$ is coherent.

Let $h'$ be a piecewise linear and strictly convex function with respect to $\Sigma'$. Since $\Sigma''$ stellarly subdivides every simplex $\sigma$ of $\Sigma$, there exist functions $h''_\sigma: \sigma \to \R$ piecewise linear and strictly convex with respect to $\Sigma''$. By an affine transformation one can achieve $h''_\sigma|_{\partial \sigma} = 0$ for every $\sigma$. Put $h'' = h' + \epsilon h''_\sigma$ on $\sigma$. Then $h''$ is piecewise linear with respect to $\Sigma''$ and strictly convex across its facets for a sufficiently small positive $\epsilon$.
\end{proof}

\begin{proof}[Outline of the proof of Theorem \ref{thm:MW}]
Applying Lemma \ref{lem:SSubd} to an arbitrary simplicial cobordism $\Sigma$ between $T$ and $T'$, we get a coherent simplicial cobordism $\Sigma''$. Since $\Sigma''|_{P \times {1}}$ and $\Sigma''|_{P \times {0}}$ are stellar subdivisions of $T$ and $T'$ respectively, this yields a sequence of stellar moves connecting $T$ and $T'$.
\end{proof}

\subsection{Realizing interior stellar moves by bistellar moves.}
In order to simplify our task in the next section, we show that instead of interior stellar moves one can use \emph{bistellar} or \emph{Pachner moves} and continuous displacements of the vertices of the triangulation.

\begin{figure}[ht]
\includegraphics[width=0.5\textwidth]{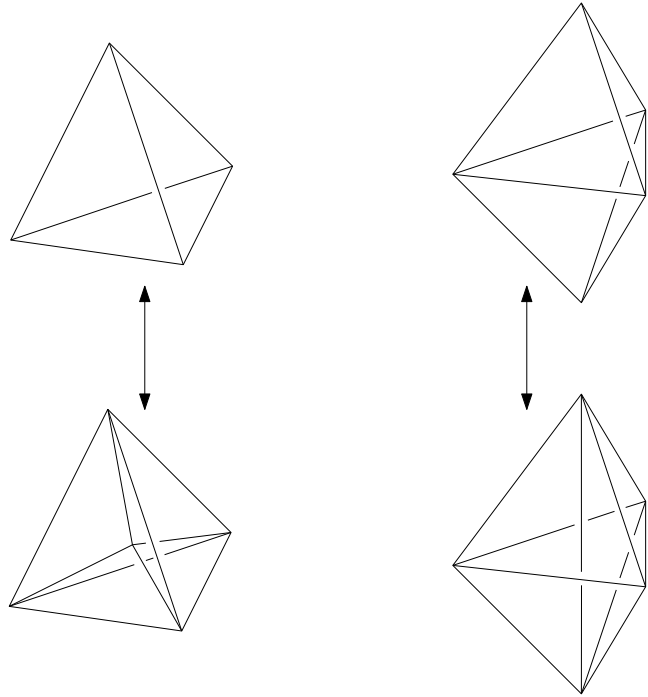}
\caption{The $1 \leftrightarrow 4$ and $2 \leftrightarrow 3$ Pachner moves.}
\label{fig:Pachner}
\end{figure}

\begin{defi}
Let $T$ be a triangulation of a subset of $\R^3$.
\begin{itemize}
\item Let $\sigma$ be a $3$-dimensional simplex of $T$. A $1 \to 4$ Pachner move replaces $\sigma$ by four smaller simplices sharing a vertex that is an interior point of $\sigma$.
\item Let $\sigma$ and $\tau$ be two $3$-simplices of $T$ such that the union $\sigma \cup \tau$ is a strictly convex bipyramid. A $2 \to 3$ Pachner move replaces $\sigma$ and $\tau$ by three simplices sharing the edge that joins the opposite vertices of $\sigma$ and $\tau$.
\item A $3 \to 2$ Pachner move is the inverse to a $2 \to 3$ Pachner move.
\item A $4 \to 1$ Pachner move is the inverse to a $1 \to 4$ Pachner move.
\end{itemize}
\end{defi}
The Pachner moves are depicted on Figure \ref{fig:Pachner}.

\begin{lemma} \label{lm:connect-pachner}
Any two triangulations of a convex polyhedron $P$ can be connected by a sequence of Pachner moves, boundary stellar moves and continuous displacements of the interior vertices.
\end{lemma}

\begin{proof}
Due to Theorem \ref{thm:MW} it suffices to show that every interior stellar move can be realized as a sequence of Pachner moves and vertex displacements. Since Pachner moves are invertible, we realize only interior starrings.

The starring in a $3$-simplex is a $1 \to 4$ Pachner move.

Consider the starring in a triangle, Figure \ref{fig:StarWeld}, middle. Denote the vertices of the triangle to be starred by $1,2,3$, and the two remaining vertices by $a$ and $b$. Perform a $1 \to 4$ move on the tetrahedron $a123$ and denote the new vertex by $p$. Then perform a $2 \to 3$ move on the tetrahedra $p123$ and $b123$. Finally move the vertex $p$ so that it lies in the triangle $123$.

To realize a starring of an edge, we also first perform a sequence of Pachner moves to obtain a triangulation combinatorially equivalent to the starring, and then move the new vertex. Denote by $a$ and $b$ the vertices of the edge to be starred, and denote the vertices in the link of the edge $ab$ by $1,2,\ldots,n$ in the cyclic order. Perform a $1 \to 4$ move on the tetrahedron $ab1n$. The new vertex $p$ should be chosen so that the plane of the triangle $abp$ does not pass through any other vertex. Let $(k,k+1)$ be the edge intersected by this plane. Perform a $2 \to 3$ move on the tetrahedra $ab1p$ and $ab12$, then a $2 \to 3$ move on the tetrahedra $ab2p$ and $ab23$, and so on. This sequence finishes with a $2 \to 3$ move on $ab(k-1)p$ and $ab(k-1)k$. After that apply a similar sequence of $2 \to 3$ moves on the other side starting with the tetrahedra $abnp$ and $abn(n-1)$ and finishing with $ab(k+2)p$ and $ab(k+2)(k+1)$. Finally perform a $3 \to 2$ move over the tetrahedra $abpk$, $abp(k+1)$ and $abk(k+1)$. It remains to move the vertex $p$ so that it lies on the edge $ab$.
\end{proof}

\section{The effect of the elementary moves on the signature of $M_T$}
\label{sec:Pach}
In this Section we realize the second point of the plan outlined in Section \ref{subsec:plan}. Namely, we show that if Theorem \ref{tm:withvertices} holds for some triangulation $T$, then it holds for a triangulation $T'$ that is obtained from $T$ by an elementary move. An elementary move is either a Pachner move or a boundary stellar move or a continuous displacement of the interior vertices of $T$.
\subsection{The rank of the matrix $M_T$.}
Here we prove a part of Theorem \ref{tm:withvertices}:

\begin{lemma}
\label{lem:corank}
The corank of the matrix $M_T$ equals $3m+k$, where $m$ is the number of the interior vertices and $k$ is the number of flat boundary vertices in the triangulation $T$:
$$
\dim \ker M_T = 3m+k.
$$
\end{lemma}
\begin{proof}
If $m > 0$ or $k > 0$, then it is easy to find a whole bunch of vectors in the kernel of $M_T$. Any continuous displacement of the interior vertices of $T$ changes the lengths of the interior edges, but doesn't change the angles around them, which stay equal to $2\pi$. Similarly, moving a flat boundary vertex in the direction orthogonal to its ambient face doesn't change any of the angles $\omega_i$. It does change the lengths of the boundary edges incident to this vertex, but only in the second order. It follows that the variations of interior edge lengths induced by the orthogonal displacement of a flat boundary vertex belong to the kernel of $M_T$.

Being formal, let $Q: \cV(T) \to \R^3$ be an assignment to every vertex $p_i$ of $T$ of a vector $q_i$ such that
\begin{enumerate}
 \item $q_i = 0$ if $p_i$ is a non-flat boundary vertex of $T$;
 \item $q_i \perp F_i$ if $p_i$ is a flat boundary vertex lying in the face $F_i$ of $P$.
\end{enumerate}
For every edge $ij$ of $T$ put
$$
\ell^Q_{ij} = \frac{p_i - p_j}{\|p_i - p_j\|} \cdot (q_i - q_j).
$$
It is easy to see that this formula gives the infinitesimal change of $\ell_{ij}$ that results from the infinitesimal displacements of the vertices $p_i, p_j$ by the vectors $q_i, q_j$. By the previous paragraph, $\ell^Q_{ij} \in \ker M_T$.

Let us show that the span of the vectors $\ell^Q$ has dimension $3m+k$. The correspondence between $Q$ and $\ell^Q$ is linear, and the space of assignments $Q$ with properties (1) and (2) has dimension $3m+k$, so it suffices to show that $\ell^Q = 0$ implies $Q=0$. Indeed, $\ell^Q = 0$ means that $Q$ is an infinitesimal flex of the $1$-skeleton of $T$, see Section \ref{subsec:def}. But $T$ is infinitesimally rigid, since every simplex is. Thus $\ell^Q = 0$ implies that $Q$ is trivial. Since $q_i = 0$ on the vertices of $P$, we have $Q = 0$.

It remains to show that any vector $\dot\ell \in \ker M_T$ has the form $\ell^Q$ for some $Q$. Let $p_1 p_2 p_3$ be a triangle of $T$. Choose $q_1, q_2,$ and $q_3$ arbitrarily. Let $p_4$ be a vertex such that there is a simplex $p_1p_2p_3p_4$ in $T$. The values of $\dot\ell_{i4}, i = 1,2,3$ determine uniquely a vector $q_4$ such that $\dot\ell_{i4} = \ell^Q_{i4}$ for $i = 1,2,3$. If $ij$ is a boundary edge of $T$, then we put $\dot\ell_{ij} = 0$. Similarly, we define $q_5$ for the vertex $p_5$ of a simplex that shares a face with $p_1p_2p_3p_4$. Proceeding in this manner, we can assign a vector $q_i$ to every vertex $p_i$, if we show that this is well-defined (we extend our assignment along paths in the dual graph of $T$, and it needs to be shown that the extension does not depend on the choice of a path). It is not hard to see that this is ensured by the property $M_T \dot\ell = 0$. Thus we have constructed an assignment $Q: \cV(T) \to \R^3$ such that $\dot\ell = \ell^Q$. Since $\dot\ell_{ij} = 0$ for every boundary edge $ij$ of $T$, the vectors $(q_i)|_{p_i \in \partial P}$ define an infinitesimal flex of the boundary of $P$. Due to Theorem \ref{thm:DehnGen}, $Q$ satisfies properties (1) and (2) above, after subtracting an infinitesimal motion. Thus the kernel of $M_T$ consists of the vectors of the form $\ell^Q$.
\end{proof}

\begin{cor}
Let $T$ be a triangulation of a convex polyhedron $P$. Consider a continuous displacement of the vertices of $T$ such that no simplex of the triangulation degenerates, the underlying space of $T$ remains a convex polyhedron, all flat boundary vertices remain flat, and non-flat remain non-flat. Then the signature of the matrix $M_T$ does not change during this deformation.
\end{cor}
\begin{proof}
Due to Lemma \ref{lem:corank}, the rank of $M_T$ does not change during the deformation. Hence no eigenvalue changes its sign.
\end{proof}

\subsection{The effect of the Pachner moves.}

\begin{lemma}\label{lm:23}
Let $P$ be a convex polyhedron, and let $T$ and $T'$ be two triangulations
of $P$ such that $T'$ is obtained from $T$ by a $2 \to 3$ Pachner move.
Then the statement of Theorem \ref{tm:withvertices} applies to $T$ if
and only if it applies to~$T'$.
\end{lemma}

\begin{proof}
Since triangulations $T$ and $T'$ have the same number of interior and flat boundary vertices, the matrices $M_T$ and $M_{T'}$ have the same corank by Lemma \ref{lem:corank}. It remains to show that $M_T$ and $M_{T'}$ have the same number of negative eigenvalues.

Matrices $M_T$ and $M_{T'}$ define symmetric bilinear forms (that are denoted by the same letters) on the spaces $\R^{\cE_{\inn}(T)}$ and $\R^{\cE_{\inn}(T')}$, respectively. Here $\cE_{\inn}(T)$ denotes the set of interior edges of the triangulation $T$. Note that $\cE_{\inn}(T') = \cE_{\inn}(T) \cup \{e_0\}$, where $e_0$ is the vertical edge on the lower right of Figure \ref{fig:Pachner}. Extend $M_T$ to a symmetric bilinear form on $\R^{\cE_{\inn}(T')}$ by augmenting the matrix $M_T$ with a zero row and a zero column, and put
$$
\Phi = M_{T'} - M_T.
$$
By Definition \ref{dfn:M_T}, we have
$$
\Phi = \left( \frac{\partial(\omega'_i - \omega_i)}{\partial \ell_j} \right)_{i,j \in \cE_{\inn}(T')},
$$
where we put $\frac{\partial\omega_0}{\partial \ell_j} = 0$ for all $j$.

Denote those edges on the upper right of Figure \ref{fig:Pachner} that are interior edges of $T$ by $e_1, \ldots, e_s$. Note that $\omega_i = \omega'_i$ as functions of the edge lengths for all $i \notin \{0,\ldots,s\}$. Thus, the matrix $\Phi$ reduces to an $(s+1) \times (s+1)$ matrix with rows corresponding to the edges $e_0, \ldots, e_s$.

We claim that the matrix $\Phi$ is positive semidefinite of rank $1$. To construct a vector in the kernel of $\Phi$, note that during any continuous deformation of the bipyramid on Figure \ref{fig:Pachner} we have $\omega_i = \omega'_i$ as functions of edge lengths for $i = 1,\ldots, s$, while $\omega'_0$ is identically $2\pi$. Thus if we choose $\dot\ell_1, \ldots, \dot\ell_s$ arbitrarily and define $\dot\ell_0$ as the infinitesimal change of the length of $e_0$ under the corresponding infinitesimal deformation of the bipyramid, then we have $\Phi \dot\ell = 0$. Therefore $\rk \Phi \le 1$. The infinitesimal rigidity of the bypiramid implies $\frac{\partial \omega'_0}{\partial \ell_0} \ne 0$, thus we have
$$
\rk \Phi = 1.
$$
Since the space of convex bypiramids is connected, it suffices to prove the positive semidefiniteness of $\Phi$ in some special case. In the case when all edges of the bypiramid have equal length one can easily see that $\frac{\partial \omega'_0}{\partial \ell_0} > 0$, which implies the positivity of the unique eigenvalue of $\Phi$.

The equation
$$
\rk M_{T'} = \rk M_T + 1 = \rk M_T + \rk \Phi
$$
implies that $\ker M_T$ and $\ker \Phi$ intersect transversally and $\ker M_{T'} = \ker M_T \cap \ker \Phi$. Therefore
$$
\rk (M_T + t\Phi) = \rk M_T + 1
$$
for all $t \ne 0$. The Courant minimax principle \cite[Chapter I, \S 4]{CH53} implies that the eigenvalues of $M_T + \epsilon \Phi$ are larger than or equal to the corresponding eigenvalues of $M_T$. It follows that when $M_T$ is deformed into $M_{T'}$ via $\{M_T + t\Phi\}_{t \in [0,1]}$, exactly one of the zero eigenvalues of $M_T$ becomes positive, and all of the non-zero eigenvalues preserve their sign. Thus $M_{T'}$ has the same number of negative eigenvalues as $M_T$ and the lemma is proved.
\end{proof}

\begin{lemma}\label{lm:14}
Let $P$ be a convex polyhedron, and let $T$ and $T'$ be two triangulations
of $P$ such that $T'$ is obtained from $T$ by a $1 \to 4$ Pachner move.
Then the statement of Theorem \ref{tm:withvertices} applies to $T$ if
and only if it applies to~$T'$.
\end{lemma}

\begin{proof}
The same arguments as in the proof of Lemma \ref{lm:23} work. The triangulation $T'$ has one interior vertex more that the triangulation $T$ and four interior edges more than $T$. Due to Lemma \ref{lem:corank}, we have
$$
\rk M_{T'} = \rk M_T + 1,
$$
and we have to prove that $M_{T'}$ has the same number of positive eigenvalues as $M_T$ and one negative eigenvalue more. For this it suffices to show that the quadratic form
$$
\Phi = M_{T'} - M_T
$$
is negative semidefinite of rank $1$. In the same way as in the proof of Lemma \ref{lm:23}, one shows that $\rk \Phi \le 1$. After that, it suffices to show that the restriction of $\Phi$ to the space spanned by the variations of lengths of the four interior edges on the lower left of Figure \ref{fig:Pachner} is non-trivial and negative semidefinite. The non-triviality follows from the infinitesimal rigidity of the simplex, and it suffices to check the negative semidefiniteness in some convenient special case.
\end{proof}

\subsection{The effect of the boundary stellar moves.}

\begin{lemma}
Let $P$ be a convex polyhedron, and let $T$ and $T'$ be two triangulations
of $P$ such that $T'$ is obtained from $T$ by the starring of a boundary $2$-simplex.
Then the statement of Theorem \ref{tm:withvertices} applies to $T$ if
and only if it applies to~$T'$.
\end{lemma}
\begin{proof}
We have
$$
\rk M_{T'} = \rk M_T
$$
and need to show that $M_{T'}$ has the same signature as $M_T$. This is true because in fact 
$M_{T'} = M_T$ -- more precisely, $M_{T'}$ is obtained from $M_T$ by adding a column and
a row, each with all elements equal to zero. 
This can be shown using the explicit formulas for $\frac{\partial \omega_i}{\partial \ell_j}$ 
and $\frac{\partial \omega'_i}{\partial \ell_j}$ from \cite[Section 3.1]{bobenko-izmestiev} and 
\cite{Kor00}.
\end{proof}

\begin{lemma}
Let $P$ be a convex polyhedron, and let $T$ and $T'$ be two triangulations
of $P$ such that $T'$ is obtained from $T$ by the starring of a boundary $1$-simplex.
Then the statement of Theorem \ref{tm:withvertices} applies to $T$ if
and only if it applies to~$T'$.
\end{lemma}
\begin{proof}
The strategy is the same as in the proofs of Lemmas \ref{lm:23} and \ref{lm:14}. Put
$$
\Phi = M_{T'} - M_T
$$
and note that by Lemma \ref{lem:corank}
\begin{equation}
\label{eqn:rk_d1}
\rk M_{T'} = \rk M_T + i,
\end{equation}
where $i$ is one less than the number of simplices incident to the starred edge (e.g. in the right column of Figure~\ref{fig:dStarWeld}, $i = 2$). As in the proof of Lemma \ref{lm:23}, one shows that $\rk \Phi \le i$. Then (\ref{eqn:rk_d1}) implies
$$
\rk \Phi = i.
$$
Since we aim to show that $M_{T'}$ has the same number of negative eigenvalues as $M_T$, it suffices to show that $\Phi$ is positively semidefinite.

Let $\Psi$ be the $i \times i$ principal minor of $\Phi$ formed by the rows and columns that correspond to the interior edges of the triangulation on the lower right of Figure \ref{fig:dStarWeld}. We claim that $\Psi$ is positively definite, which implies the nonnegativity of $\Phi$. The proof is by continuity argument as in Lemma \ref{lm:23}. In order to prove the non-degeneracy of $\Psi$, it suffices to show that the framework of the boundary edges on the lower right of Figure \ref{fig:dStarWeld} is infinitesimally rigid. Note that the framework on the upper right of Figure \ref{fig:dStarWeld} is infinitesimally rigid, since it is formed by skeleta of $3$-simplices that are rigid. This implies the infinitesimal rigidity of the boundary framework on the lower right (an easy exercise in applying the definition of an infinitesimal flex). Now consider a deformation of the triangulation on the upper right that makes the underlying polyhedron
convex. This deformation can be extended to a deformation of the triangulation on the lower right. Since the matrix $\Psi$ remains non-degenerate during the deformation, its signature is preserved. After the polyhedron is made convex, push the starring vertex off the starred edge so that the vertices of the triangulation are in the convex position. This also preserves the signature of $\Psi$. In the final position $\Psi$ is positive due to Theorem \ref{thm:StarPos}.
\end{proof}

\section{Investigating $M_T$ for a special triangulation $T$.}
\label{sec:Special}

Let $P$ be a convex polyhedron. Let $S$ be a triangulation of $\partial P$ with $\V(S) = \V(P)$, 
and let $p$ be a vertex of $P$. Consider the triangulation $T$ that consists of simplices with 
a common vertex $p$ and opposite faces the triangles of $S$ disjoint with $p$.
\begin{thm}
\label{thm:StarPos}
The matrix $M_T$ is positive definite.
\end{thm}
\begin{proof}
This is proved in \cite{star}. Formally, it is a special case of \cite[Theorem 1.5]{star} that 
claims that $M_T$ is positive if $P$ is weakly convex and star-shaped with respect to the vertex 
$p$. The proof uses the positivity of the corresponding matrix for convex caps 
\cite[Lemma 6, Theorem 5]{izmestiev} and the projective invariance of the infinitesimal rigidity \cite[Section 5]{star}.
\end{proof}

Theorem \ref{thm:StarPos} accomplishes the plan outlined in Section \ref{subsec:plan}. 
Theorem \ref{tm:withvertices} is proved, and therewith Theorem~\ref{tm:main}.

\appendix

\section{A polyhedron which is not weakly codecomposable}

We describe in this appendix a simple example of a polyhedron which is not
weakly codecomposable in the sense of Definition \ref{dfn:codec}. 

\begin{defi}
Let $\theta\in (-2\pi/3,2\pi/3)$. The twisted octahedron $Oct_\theta$ is the polyhedron
with vertices $A,B,C,A',B',C'$ of coordinates respectively $(1,0,1)$,
$(\cos(2\pi/3),\sin(2\pi/3),1)$, $(\cos(4\pi/3),\sin(4\pi/3),1)$,
$(\cos(-\pi+\theta), \sin(-\pi+\theta), -1)$, 
$(\cos(-\pi/3+\theta), \sin(-\pi/3+\theta), -1)$, 
$(\cos(\pi/3+\theta), \sin(\pi/3+\theta), -1)$.
The edges are the segments joining $A$ to $B'$ and $C'$, $B$ to 
$A'$ and $C'$, $C$ to $A'$ and $B'$, and the faces are the triangles 
$(ABC), (A'B'C'),(AB'C'),(A'BC'),(A'B'C),(ABC'),(AB'C),(A'BC)$.
\end{defi}

Note that $Oct_{\pm\frac{\pi}{2}}$ is a Wunderlich's twisted octahedron, see Figure \ref{fig:IcoOcta}, right.

\begin{prop}
$Oct_\theta$ is embedded for all $\theta\in (-2\pi/3,2\pi/3)$.
\end{prop}

For $\theta\in (-2\pi/3,2\pi/3)$ we call $A_t(\theta)$ the area of the
intersection of $Oct_\theta)$ with the horizontal plane $\{ z=t\}$.

\begin{prop}
$\lim_{\theta\rightarrow 2\pi/3} A_0(\theta)=0$. 
\end{prop}

Let $K$ be a large enough convex polygon in the plane $Oxy$ 
(it suffices to require that the interior of $K$ contains the disk $x^2+y^2 \le 1$). 
Consider the polyhedron $P_\theta = \conv(A,B,C,A',B',C',K) \setminus Oct_\theta$ homeomorphic to a solid torus.

\begin{lemma}
For $\theta$ close enough to $2\pi/3$, $P_\theta$ is not weakly 
codecomposable.
\end{lemma}

\begin{proof}
Suppose that $P_\theta$ is weakly codecomposable. Then there exists a 
convex polyhedron $Q_\theta\supset P_\theta$ such that 
 $Q_\theta\setminus P_\theta$ can be triangulated without interior vertex.
Let $S_1,\cdots, S_n$ be the simplices in this triangulation which intersect
$Oct_\theta\cap (Oxy)$. For each $i\in \{ 1,\cdots, n\}$, let $a_i(t)$ be the 
area of the intersection of $S_i$ with the horizontal plane $\{ z=t\}$.

Each of the $S_i$ can have either:
\begin{itemize}
\item two vertices with $z\geq 1$ and two vertices with $z\leq -1$. Then 
the restriction of $a_i$ to $(-1,1)$ is a concave
quadratic function, so that $2a_i(0)\geq a_i(-1)+a_i(1)$.
\item one vertex with $z\geq 1$ and three vertices with $z\leq -1$. Then
$a_i$ is of the form $a_i(t) = c_i(t+b_i)^2$ with $b_i \ge 1$. It easily implies that 
$4a_i(0)\geq a_i(-1)+a_i(1)$.
\item one vertex with $z\leq -1$ and three vertices with $z\geq 1$. The
same argument then shows the same result.
\end{itemize}

So $4a_i(0)\geq a_i(-1)+a_i(1)$ for all $i$ and the union of the $S_i$
contains $Oct_\theta$. It follows that $4A_0(\theta)\geq A_{-1}(\theta)
+A_1(\theta)$. But $A_1(\theta)$ and $A_{-1}(\theta)$ are equal to the
area of an equilateral triangle of fixed side length, while $A_0(\theta)$
goes to $0$ as $\theta\rightarrow 2\pi/3$, this is a contradiction.
So $P_\theta$ is not weakly codecomposable for $\theta$ close enough
to $2\pi/3$.
\end{proof}


\providecommand{\bysame}{\leavevmode\hbox to3em{\hrulefill}\thinspace}
\providecommand{\MR}{\relax\ifhmode\unskip\space\fi MR }
\providecommand{\MRhref}[2]{%
  \href{http://www.ams.org/mathscinet-getitem?mr=#1}{#2}
}
\providecommand{\href}[2]{#2}

\end{document}

%% file: macros.tex
\newtheorem{cor}{Corollary}[section]
\newtheorem{thm}[cor]{Theorem}
\newtheorem{conj}[cor]{Conjecture}
\newtheorem{prop}[cor]{Proposition}
\newtheorem{lemma}[cor]{Lemma}
\theoremstyle{definition}
\newtheorem{defi}[cor]{Definition}
\theoremstyle{remark}
\newtheorem{remark}[cor]{Remark}
\newtheorem{example}[cor]{Example}

\newcommand{\cD}{{\mathcal D}}
\newcommand{\cM}{{\mathcal M}}
\newcommand{\cS}{{\mathcal S}}
\newcommand{\cT}{{\mathcal T}}
\newcommand{\cML}{{\mathcal M\mathcal L}}
\newcommand{\cGH}{{\mathcal G\mathcal H}}
\newcommand{\C}{{\mathbb C}}
\newcommand{\N}{{\mathbb N}}
\newcommand{\R}{{\mathbb R}}
\newcommand{\Z}{{\mathbb Z}}
\newcommand{\Kt}{\tilde{K}}
\newcommand{\Mt}{\tilde{M}}
\newcommand{\dr}{{\partial}}
\newcommand{\kappab}{\overline{\kappa}}
\newcommand{\pib}{\overline{\pi}}
\newcommand{\Sigmab}{\overline{\Sigma}}
\newcommand{\gd}{\dot{g}}
\newcommand{\ld}{\dot{l}}
\newcommand{\thetad}{\dot{\theta}}
\newcommand{\diff}{\mbox{Diff}}
\newcommand{\dev}{\mbox{dev}}
\newcommand{\devb}{\overline{\mbox{dev}}}
\newcommand{\devt}{\tilde{\mbox{dev}}}
\newcommand{\vol}{\mbox{Vol}}
\newcommand{\hess}{\mbox{Hess}}
\newcommand{\db}{\overline{\partial}}
\newcommand{\Sigmat}{\tilde{\Sigma}}

\newcommand{\cunc}{{\mathcal C}^\infty_c}
\newcommand{\cun}{{\mathcal C}^\infty}
\newcommand{\dd}{d_D}
\newcommand{\dmin}{d_{\mathrm{min}}}
\newcommand{\dmax}{d_{\mathrm{max}}}
\newcommand{\Dom}{\mathrm{Dom}}
\newcommand{\dn}{d_\nabla}
\newcommand{\ded}{\delta_D}
\newcommand{\delmin}{\delta_{\mathrm{min}}}
\newcommand{\delmax}{\delta_{\mathrm{max}}}
\newcommand{\hmin}{H_{\mathrm{min}}}
\newcommand{\maxi}{\mathrm{max}}
\newcommand{\oL}{\overline{L}}
\newcommand{\oP}{{\overline{P}}}
\newcommand{\Ran}{\mathrm{Ran}}
\newcommand{\tgamma}{\tilde{\gamma}}
\newcommand{\cotan}{\mbox{cotan}}
\newcommand{\lambdat}{\tilde\lambda}
\newcommand{\St}{\tilde S}

\newcommand{\II}{I\hspace{-0.1cm}I}
\newcommand{\III}{I\hspace{-0.1cm}I\hspace{-0.1cm}I}
\newcommand{\note}[1]{\marginpar{\tiny #1}}

\newcommand{\cV}{{\mathcal V}}
\newcommand{\cE}{{\mathcal E}}
\newcommand{\inn}{\mathrm{int}}
\newcommand{\rk}{\mathrm{rank\,}}
\newcommand{\st}{\mathrm{st\,}}
\newcommand{\V}{\mathrm{Vert}}
\newcommand{\conv}{\mathrm{conv\,}}
\newcommand{\pr}{\mathrm{pr}}

%% file: weakly.bbl
\providecommand{\bysame}{\leavevmode\hbox to3em{\hrulefill}\thinspace}
\providecommand{\MR}{\relax\ifhmode\unskip\space\fi MR }
\providecommand{\MRhref}[2]{%
  \href{http://www.ams.org/mathscinet-getitem?mr=#1}{#2}
}
\providecommand{\href}[2]{#2}
\begin{thebibliography}{10}

\bibitem{MW07}
\emph{Jessen's orthogonal icosahedron}, {\tt
  http://mathworld.wolfram.com/JessensOrthogonalIcosahedron.html}.

\bibitem{aichholzer}
Oswin Aichholzer, Lyuba~S. Alboul, and Ferran Hurtado, \emph{On flips in
  polyhedral surfaces}, Internat. J. Found. Comput. Sci. \textbf{13} (2002),
  no.~2, 303--311, Volume and surface triangulations.

\bibitem{alex}
A.~D. Alexandrov, \emph{Convex polyhedra}, Springer Monographs in Mathematics,
  Springer-Verlag, Berlin, 2005.

\bibitem{bobenko-izmestiev}
A.~Bobenko and I.~Izmestiev, \emph{Alexandrov's theorem, weighted {D}elaunay
  triangulations, and mixed volumes}, {\tt arXiv:math.DG/0609447}, to appear in
  Ann. Inst. Fourier, 2007.

\bibitem{cauchy}
Augustin~Louis Cauchy, \emph{Sur les polygones et poly\`edres, second
  m\'emoire}, Journal de l'Ecole Polytechnique \textbf{19} (1813), 87--98.

\bibitem{vienna}
Robert Connelly and Jean-Marc Schlenker, \emph{On the infinitesimal rigidity of
  weakly convex polyhedra}, {\tt arXiv:math.DG/0606681}, 2006.

\bibitem{CH53}
R.~Courant and D.~Hilbert, \emph{Methods of mathematical physics. {V}ol. {I}},
  Interscience Publishers, Inc., New York, N.Y., 1953.

\bibitem{dehn}
Max Dehn, \emph{{\"U}ber die {S}tarrheit konvexer {P}olyeder.}, Math. Ann.
  \textbf{77} (1916), 466--473.

\bibitem{glaser}
L.~C. Glaser, \emph{Geometric combinatorial topology}, Van Nostrand Reinhold,
  New York, 1970.

\bibitem{Glu75}
Herman Gluck, \emph{Almost all simply connected closed surfaces are rigid},
  Geometric topology (Proc. Conf., Park City, Utah, 1974), Springer, Berlin,
  1975, pp.~225--239. Lecture Notes in Math., Vol. 438.

\bibitem{izmestiev}
Ivan Izmestiev, \emph{A variational proof of {Alexandrov's} convex cap
  theorem}, {\tt arXiv:math.DG/0703169}, 2007.

\bibitem{Jes67}
B{\o}rge Jessen, \emph{Orthogonal icosahedra}, Nordisk Mat. Tidskr \textbf{15}
  (1967), 90--96.

\bibitem{Kor00}
I.~G. Korepanov, \emph{A formula with volumes of five tetrahedra and discrete
  curvature}, {\tt arXiv:nlin/0003001}, 2000.

\bibitem{legendre}
A.~M. Legendre, \emph{El\'ements de g\'eom\'etrie}, Paris, 1793 (an II),
  Premi\`ere \'edition, note XII, pp. 321--334.

\bibitem{Lick99}
W.~B.~R. Lickorish, \emph{Simplicial moves on complexes and manifolds},
  Proceedings of the Kirbyfest (Berkeley, CA, 1998), Geom. Topol. Monogr.,
  vol.~2, Geom. Topol. Publ., Coventry, 1999, pp.~299--320 (electronic).

\bibitem{morelli}
Robert Morelli, \emph{The birational geometry of toric varieties}, J. Algebraic
  Geom. \textbf{5} (1996), no.~4, 751--782.

\bibitem{regge}
T~Regge, \emph{General relativity without coordinates}, Nuovo Cimento
  \textbf{19} (1961), 558--571.

\bibitem{San05}
Francisco Santos, \emph{Non-connected toric {H}ilbert schemes}, Math. Ann.
  \textbf{332} (2005), no.~3, 645--665.

\bibitem{San06}
\bysame, \emph{Geometric bistellar flips: the setting, the context and a
  construction}, International Congress of Mathematicians. Vol. III, Eur. Math.
  Soc., Z\"urich, 2006, pp.~931--962.

\bibitem{rcnp}
Jean-Marc Schlenker, \emph{A rigidity criterion for non-convex polyhedra},
  Discrete Comput. Geom. \textbf{33} (2005), no.~2, 207--221.

\bibitem{star}
\bysame, \emph{On weakly convex star-shaped polyhedra}, {\tt arXiv:0704.2901},
  2007.

\bibitem{wlodarczyk}
Jaros{\l}aw W{\l}odarczyk, \emph{Decomposition of birational toric maps in
  blow-ups \& blow-downs}, Trans. Amer. Math. Soc. \textbf{349} (1997), no.~1,
  373--411. \MR{MR1370654 (97d:14021)}

\bibitem{Wun65}
W.~Wunderlich, \emph{Starre, kippende, wackelige und bewegliche
  {A}chtfl\"ache}, Elem. Math. \textbf{20} (1965), 25--32.

\end{thebibliography}
